\documentclass{IEEEtran4PSCC}
\ifCLASSINFOpdf
   \usepackage[pdftex]{graphicx}
\else
   \usepackage[dvips]{graphicx}
\fi
\usepackage{multirow}

\usepackage[cmex10]{amsmath}
 \usepackage[caption=false,font=footnotesize]{subfig}
\usepackage{times}
\usepackage{graphicx} 
\usepackage{tikz}
\usepackage{amsmath}
\usepackage{amssymb}

\usepackage{amsthm}

\usepackage{draftfigure} 

\DeclareMathOperator*{\argmin}{arg\,min}
\DeclareMathOperator*{\minimize}{min}
\DeclareMathOperator*{\subjectto}{s.t.}

\DeclareMathOperator*{\diag}{diag}

\newcommand{\dd}{\mathsf{d}}

\makeatletter
\let\old@ps@headings\ps@headings
\let\old@ps@IEEEtitlepagestyle\ps@IEEEtitlepagestyle
\def\psccfooter#1{%
    \def\ps@headings{%
        \old@ps@headings%
        \def\@oddfoot{\strut\hfill#1\hfill\strut}%
        \def\@evenfoot{\strut\hfill#1\hfill\strut}%
    }%
    \def\ps@IEEEtitlepagestyle{%
        \old@ps@IEEEtitlepagestyle%
        \def\@oddfoot{\strut\hfill#1\hfill\strut}%
        \def\@evenfoot{\strut\hfill#1\hfill\strut}%
    }%
    \ps@headings%
}
\makeatother

\psccfooter{%
        \parbox{\textwidth}{\hrulefill \\ \small{24th Power Systems Computation Conference} \hfill \begin{minipage}{0.2\textwidth}\centering \vspace*{4pt} \includegraphics[scale=0.06]{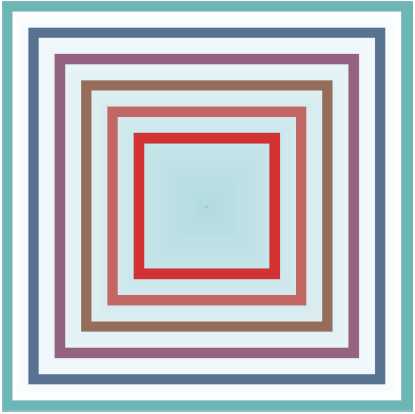}\\\small{PSCC 2026} \end{minipage} \hfill \small{Limassol, Cyprus --- June 8-12, 2026}}%
}

\begin{document}
%
\title{Self-Supervised Learning of Parametric Approximation for Security-Constrained DC-OPF}

\author{
\IEEEauthorblockN{Anderson Anrrango}
\IEEEauthorblockA{Department of Electrical Energy\\
Escuela Politécnica Nacional\\
Quito, Ecuador
}
\and
\IEEEauthorblockN{André H. Quisaguano, Gonzalo E. Constante, Can Li}
\IEEEauthorblockA{Davidson School of Chemical Engineering \\
Purdue University\\
West Lafayette, IN, USA
}
}


\maketitle

\begin{abstract}
This paper introduces a self-supervised learning framework for approximating the Security-Constrained DC Optimal Power Flow (SC-DCOPF) problem using a parametric linear model. The approach preserves the physical structure of the DC-OPF while incorporating demand-dependent tunable parameters that scale transmission line limits. These parameters are predicted via a Graph Neural Network and optimized through differentiable layers, enabling direct training from contingency costs without requiring labeled data. The framework integrates pre- and post-contingency optimization layers into an implicit loss function. Numerical experiments on benchmark systems demonstrate that the proposed method achieves high dispatch accuracy, low cost approximation error, and strong data efficiency, outperforming semi-supervised and end-to-end baselines. This scalable and interpretable approach offers a promising solution for real-time secure power system operations.
\end{abstract}


\thanksto{\noindent Submitted to the 24th Power Systems Computation Conference (PSCC 2026).}

\section{Introduction}
\subsection{Context and motivation}

The Security-Constrained Optimal Power Flow (SC-DCOPF) problem is a core problem in power system operation that ensures secure and economic dispatch under possible component failures \cite{Capitanescu2011}. However, solving the SC-DCOPF problem is computationally demanding due to the size of the optimization problem, whose number of variables and constraints grow rapidly as a function of the number of contingencies. To address this, recent works have explored various approximation and learning-based strategies to accelerate SC-DCOPF solutions while maintaining acceptable accuracy \cite{Pan2021,Park2025}. Nevertheless, these approaches often come at the cost of reduced interpretability or lack formal guarantees of constraint satisfaction \cite{Kotary2021,Ferrando2024}, which limits their applicability in safety-critical operational settings.

Recent advances in machine learning have enabled the development of parametric optimization models that combine physics-based formulations with data-driven parameterization \cite{Amos2017,Agrawal2019,Wilder2019,Shah2022}. These models accelerate decision-making, improve accuracy, and preserve interpretability.

\subsection{Literature review}

In the context of the optimal power flow (OPF) problem, Taheri and Molzahn \cite{Taheri2024,Taheri2025} propose tuning DC power flow parameters, such as susceptance coefficients and biases, via gradient-based methods to improve approximation accuracy relative to the AC formulation over specific operating ranges. Constante-Flores et al. \cite{Constante2025} adopt a supervised learning approach to map operating conditions to optimal tunable parameters, enhancing DC-OPF dispatch accuracy while maintaining desirable market properties. Similarly, Rosemberg and Klamkin \cite{Rosemberg2025} introduce a semi-supervised learning method to map operating conditions directly to the optimal generation dispatch of a parameterized DC-OPF problem using a differentiable optimization layer. In the aforementioned frameworks, the tunable parametric models are obtained by minimizing the Euclidean distance with respect to the AC-OPF solution. However, relying on such accuracy-based metrics can yield suboptimal outcomes, both in optimality and feasibility \cite{Elmachtoub2022}. To bridge this gap, recent works \cite{Chen2024,Chen2025} have highlighted the advantage of using more realistic loss functions in the context of the ACOPF problem that capture automatic generator responses as well as asymmetries of the prediction errors in the operation of power systems.

\subsection{Contributions}

We propose a parametric linear approximation of the SC-DCOPF problem that retains the physical structure of the DC-OPF model while introducing demand-dependent tunable parameters to scale transmission line limits and approximate pre-contingency dispatch. These parameters form an inner approximation of the DC-OPF and are learned via a self-supervised framework with differentiable optimization layers, allowing direct tuning from contingency costs without requiring precomputed SC-DCOPF solutions. The framework consists of three components: (i) a machine learning model that predicts tunable parameters from demand profiles, (ii) a differentiable optimization layer implementing the parametric DC-OPF, and (iii) an implicit loss function that evaluates the approximate pre-contingency dispatch across all plausible N-1 contingencies and computes the total cost as the training objective. The proposed method trains parameters to minimize the pre- and post-contingency costs and thereby obtain a faithful approximation of the secure dispatch.

The remainder of this paper is organized as follows: Section~\ref{sec:proposed} presents the problem setup and introduces the proposed approach. Section~\ref{sec:approach} details the proposed learning methodology as well as the computation of sensitivities. Section~\ref{sec:experiments}  illustrates the performance of the proposed approach for a number of instances of the PowerGrid Library. 
Section~\ref{sec:conclusion} concludes the paper and draws future research directions.

\section{Linear Proxy of the SC-DCOPF Problem \label{sec:proposed}}

\subsection{Preliminaries}
This subsection presents the formulation of the DC OPF problem and its corrective security-constrained counterpart. 
We consider a transmission network with $n$ buses,  $l$ transmission lines, and $g$ generating units.

The DC-OPF problem can be formulated as follows:
\begin{subequations}\label{eq:dcopf}
\begin{align}
\minimize_{\mathbf{p}} \;\; & \tfrac{1} {2}\mathbf{p}^\top \mathbf{Cp} + \mathbf{c}^\top \mathbf{p} \label{eq:dcopf_a}\\
\rm{s.t.} \;\; & \mathbf{1}^\top \mathbf{p} = \mathbf{1}^\top \mathbf{d}, \label{eq:dcopf_b}\\
& \mathbf{0} \le \mathbf{p} \le \overline{\mathbf{p}}, \label{eq:dcopf_c}\\
& \lvert \mathbf{M\bigl(Ap -d\bigr)} \rvert \le \overline{\mathbf{f}}, \label{eq:dcopf_d}
\end{align}
\end{subequations}
where $\mathbf{p}\in \mathbb{R}^g$ is the vector of active power outputs of the dispatchable generators,  
$\mathbf{d}\in \mathbb{R}^n$ is the vector of active power net demands,  
$\mathbf{A}\in \mathbb{R}^{n \times g}$ is an incidence matrix mapping generators to buses,  
$\mathbf{M}\in \mathbb{R}^{l \times n}$ is the Power Transfer Distribution Factor (PTDF) matrix,  
and $\mathbf{1}$ (resp. $\mathbf{0}$) denotes a vector of ones (resp. zeros) of appropriate dimension.  In the objective function \eqref{eq:dcopf_a}, $\mathbf{c}$ denotes the linear generation cost coefficients, while $\mathbf{C}$ is a positive diagonal matrix of quadratic cost coefficients. The vector $\overline{\mathbf{p}} \in \mathbb{R}^g$ corresponds to the maximum generation capacities, and $\overline{\mathbf{f}} \in \mathbb{R}^l$ specifies the transmission line flow limits.

The objective function \eqref{eq:dcopf_a} is the total generation cost, modeled as a convex quadratic function of generator outputs. Constraints \eqref{eq:dcopf_b} enforce the balance between total generation and total demand, \eqref{eq:dcopf_c} ensure that each generator operates within its physical lower and upper capacity bounds, and \eqref{eq:dcopf_d} limit the transmission line flows.

The corrective approach of the SC-DCOPF problem can be formulated as follows:
\begin{subequations}\label{eq:scopf}
\begin{align}
\minimize_{\mathbf{p},\mathbf{p}^k,\mathbf{s}^k} \;\;
& \tfrac{1}{2}\mathbf{p}^\top \mathbf{C}\mathbf{p} + \mathbf{c}^\top \mathbf{p} + \tfrac{\rho}{|\mathcal{K}|} \sum_{k \in \mathcal{K}}\mathbf{1}^\top \mathbf{s}^k \label{eq:scopf_a}\\
\rm{s.t.} \;\;
& \mathbf{1}^\top \mathbf{p} = \mathbf{1}^\top \mathbf{d}, \label{eq:scopf_b}\\
& \mathbf{0} \le \mathbf{p} \le \overline{\mathbf{p}}, \label{eq:scopf_c}\\
& \lvert \mathbf{M\bigl(Ap -d\bigr)} \rvert \le \overline{\mathbf{f}}, \label{eq:scopf_d}\\
& \forall k \in \mathcal{K}: \nonumber \\
& \qquad \mathbf{1}^\top \mathbf{p}^k = \mathbf{1}^\top (\mathbf{d} - \mathbf{s}^k), \label{eq:scopf_e}\\
& \qquad \mathbf{0} \le \mathbf{p}^k \le \overline{\mathbf{p}}, \label{eq:scopf_f}\\
& \qquad \lvert \mathbf{M}^k\bigl(\mathbf{A}\mathbf{p}^k-(\mathbf{d}-\mathbf{s}^k)\bigr) \rvert \le \overline{\mathbf{f}}^k, \label{eq:scopf_g}\\
& \qquad \mathbf{s}^k  \ge \mathbf{0}, \label{eq:scopf_h}\\
& \qquad -\underline{\mathbf{r}} \le \mathbf{p}^k - \mathbf{p} \le \overline{\mathbf{r}}, \label{eq:scopf_i}
\end{align}
\end{subequations}
where $\mathcal{K}$ denotes the set of contingencies indexed by $k$ and $\rho$ is a nonnegative penalty coefficient. The decision variables are the pre-contingency generation dispatch $\mathbf{p}\in\mathbb{R}^g$, the post-contingency corrective dispatch $\mathbf{p}_k\in\mathbb{R}^g$ for each contingency $k\in\mathcal{K}$, and the nonnegative load-shedding slack variables $\mathbf{s}_k\in\mathbb{R}^n$. 

The objective function \eqref{eq:scopf_a} corresponds to the total base-case generation cost plus a penalization term on total load shedding across contingencies, weighted by the parameter $\rho$. Constraints \eqref{eq:scopf_b}-\eqref{eq:scopf_d} enforce pre-contingency feasibility. For each contingency $k\in\mathcal{K}$, constraints \eqref{eq:scopf_e}-\eqref{eq:scopf_h} enforce the existence of a feasible corrective redispatch $\mathbf{p}_k$. Constraints \eqref{eq:scopf_i} ensure that the corrective redispatch is bounded by the ramping limits. 

\subsection{Motivation}

One of the main challenges in solving the SC-DCOPF problem is its large-scale nature, which arises from the large number of modeled contingencies. We propose an approach that addresses this by solving a parametric optimization model whose size does not increase with the number of contingencies. This parametric model yields an optimal dispatch that closely matches the pre-contingency solution of the full SC-DCOPF problem while ensuring feasibility with respect to all pre-contingency constraints.  

To achieve this, we parameterize the DC-OPF solution using a tunable vector $\boldsymbol{\alpha} \in \mathbb{R}^l$, which scales the transmission line capacities $\overline{\mathbf{f}}$. The resulting parametric DC-OPF problem is given by  
\begin{equation}\label{eq:p_dcopf}
\begin{split}
\mathbf{p}^\star = \argmin_{\mathbf{p}} \;\; & \tfrac{1} {2}\mathbf{p}^\top \mathbf{Cp} + \mathbf{c}^\top \mathbf{p} \\
\rm{s.t.} \;\; & \mathbf{1}^\top \mathbf{p} = \mathbf{1}^\top \mathbf{d},\\
& \mathbf{0} \le \mathbf{p} \le \overline{\mathbf{p}},\\
& \lvert \mathbf{M\bigl(Ap -d\bigr)} \rvert \le \overline{\mathbf{f}} \odot \boldsymbol{\alpha},
\end{split}
\end{equation}
where $\odot$ denotes the Hadamard (elementwise) product. The tuned dispatch $\mathbf{p}^\star$ is a parametric function of $\boldsymbol{\alpha}$. 

The tunable parameter $\boldsymbol{\alpha}$  is chosen to yield a dispatch $\mathbf{p}^\star$ that renders the same optimal pre-contingency dispatch of the original SC-DCOPF problem \eqref{eq:scopf} while satisfying the constraints of the pre-contingency operation, \eqref{eq:scopf_b}-\eqref{eq:scopf_d}. Note that different problem nodal demands $\mathbf{d}$ would render different optimal values of the tunable parameter vector $\boldsymbol{\alpha}$.

To model this dependency and leverage the underlying graph structure of the SC-DCOPF problem, we approximate the mapping between the nodal demands and the optimal tunable parameter by a graph neural network (GNN) as follows:
\begin{equation}
    \boldsymbol{\alpha}^\star(\mathbf{d}) \approx \boldsymbol{\tilde{\alpha}}(\mathbf{d}) = \rm{GNN}\left(\Xi ; \mathbf{d} \right),
\end{equation}
where $\boldsymbol{\alpha}^\star$ denotes the optimal tunable parameter and $\Xi$ denotes the learnable parameters of the GNN. 

\subsection{Proposed learning task}
We propose a novel self-supervised learning approach to find the learnable parameters of the GNN, by-passing the need to create a dataset before training the GNN. The proposed learning approach can be posed as the following bilevel optimization problem:
\begin{equation}\label{eq:bilevel}
\begin{split}
\minimize_{\boldsymbol{\tilde{\alpha}}, \Xi, \mathbf{p},\mathbf{p}^k,\mathbf{s}^k} \;\; 
& \mathbb{E}_\mathbf{d} \Biggl[ \tfrac{1}{2}\mathbf{p}^\top \mathbf{C}\mathbf{p} + \mathbf{c}^\top \mathbf{p} + \tfrac{\rho}{|\mathcal{K}|} \sum_{k \in \mathcal{K}}\mathbf{1}^\top \mathbf{s}^k \Biggr] \\
\rm{s.t.} \;\;
& \forall \mathbf{d}: \\
& \qquad \eqref{eq:scopf_e}- \eqref{eq:scopf_i},\; \forall k \in \mathcal{K} \\
& \qquad \boldsymbol{\tilde{\alpha}} = \rm{GNN}\left(\Xi; \mathbf{d} \right), \\
& \qquad \mathbf{0} \le \boldsymbol{\tilde{\alpha}} \le \mathbf{1},\\
& \qquad \begin{split}
\mathbf{p} = \argmin_{\mathbf{\tilde{p}}} \;\; & \tfrac{1} {2}\mathbf{\tilde{p}}^\top \mathbf{C\tilde{p}} + \mathbf{c}^\top \mathbf{\tilde{p}} \\
\rm{s.t.} \;\; & \mathbf{1}^\top \mathbf{\tilde{p}} = \mathbf{1}^\top \mathbf{d},\\
& \mathbf{0} \le \mathbf{\tilde{p}} \le \overline{\mathbf{p}},\\
& \lvert \mathbf{M(A\tilde{p} -d)} \rvert \le \overline{\mathbf{f}} \odot \boldsymbol{\tilde{\alpha}}. 
\end{split}
\end{split}
\end{equation}

In problem \eqref{eq:bilevel}, we enforce $\mathbf{0} \le \boldsymbol{\tilde{\alpha}} \le \mathbf{1}$ to yield a lower-level constraint set that is an inner approximation of the pre-contingency constraint set, defined by \eqref{eq:scopf_b}-\eqref{eq:scopf_d}. Hence, the solution of the lower-level problem $\mathbf{p}$ is also feasible for the pre-contingency constraints.

\section{Proposed Learning Framework \label{sec:approach}}
\begin{figure}[ht]
    \centering
    \includegraphics[width=\columnwidth, trim=5 2 10 30, clip]{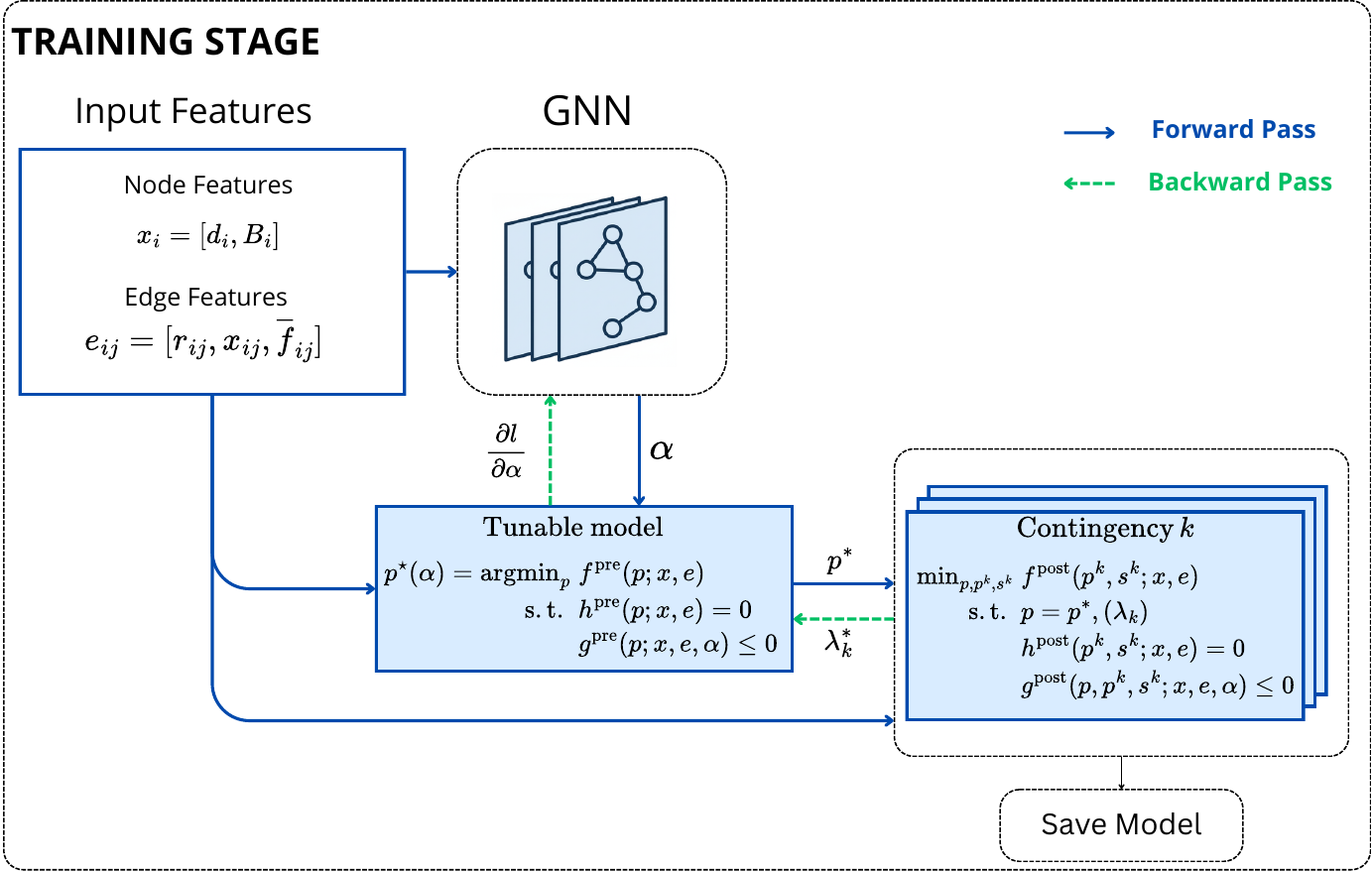}
    \caption{The proposed learning framework.}
    \label{fig:framework}
\end{figure}

The proposed framework consists of three modules, as illustrated in Fig.~\ref{fig:framework}:  
\begin{enumerate}
    \item \textit{Inference of the tunable parameter $\boldsymbol{\alpha}$}: a GNN processes the network topology and operating conditions to predict the tunable parameter $\boldsymbol{\alpha}$ that adjusts transmission line capacities.  
    \item \textit{Parametric linear approximation}: the predicted $\boldsymbol{\alpha}$ is used to solve a parametric DC-OPF model, yielding the pre-contingency generation dispatch $\mathbf{p}^\star$.  
    \item \textit{Loss evaluation}: the quality of the predicted dispatch is assessed via an implicit loss function that accounts for feasibility and optimality with respect to the original SC-DCOPF problem.  
\end{enumerate}


\subsection{Graph Neural Network}

The GNN is designed to exploit the graph structure of the power system. In this work, we adopt a Graph Attention Network v2 (GATv2) architecture, which extends the original GAT by computing attention scores through a more expressive formulation that applies the weight transformation before the attention mechanism. This design allows the network to capture richer interactions between nodes during message passing. The output layer predicts the low-dimensional parameter vector 
$\boldsymbol{\alpha} \in \mathbb{R}^l$, where $l$ is the number of transmission lines. 
This parametrization enables the optimization layer to remain independent 
of the number of contingencies, thereby addressing the scalability challenge 
of the SC-DCOPF problem.

\subsubsection{Graph encoding} 
The power system is modeled as a graph $\mathcal{G} = (\mathcal{V}, \mathcal{E})$, where $\mathcal{V}$ is the set of buses and $\mathcal{E}$ is the set of transmission lines. Each node and edge is associated with a feature vector that encodes relevant physical and operational information.

\subsubsection{Node features} 

Each node \( i \in \mathcal{V} \) represents a bus and is characterized by a feature vector containing both electrical and operational parameters. Specifically, the node features include the active power demand \( d_i \), and the corresponding network susceptance \( B_i \). These features represent the local operating conditions at each bus.

\subsubsection{Edge features} 

Each edge \( (i,j) \in \mathcal{E} \) corresponds to a transmission line, the edge features include the line resistance \( r_{ij} \), reactance \( x_{ij} \) and thermal flow limit \( \bar{f}_{ij} \). These quantities collectively describe the physical parameters that govern power transfer between buses \( i \) and \( j \).

\subsubsection{Message passing} 
In a GATv2, node embeddings are updated using attention-weighted aggregation of neighboring features.  
The update rule for node $i$ at layer $t$ is:
\begin{equation}
    h_i^{(t+1)} = \sigma\left( \sum_{j \in \mathcal{N}(i)} \gamma_{ij}^{(t)} W^{(t)} h_j^{(t)} \right),
\end{equation}
where $h_j^{(t)}$ is the embedding of neighbor node $j$, $W^{(t)}$ is a learnable weight matrix, $\sigma(\cdot)$ is a non-linear activation, and $\gamma_{ij}^{(t)}$ is the attention coefficient:
\begin{equation}
    \gamma_{ij}^{(t)} = 
    \frac{\exp\left( \text{LeakyReLU}\left( a^\top W^{(t)} [h_i^{(t)} \, || \, h_j^{(t)}] \right) \right)}
         {\sum_{k \in \mathcal{N}(i)} \exp\left( \text{LeakyReLU}\left( a^\top W^{(t)} [h_i^{(t)} \, || \, h_k^{(t)}] \right) \right)},
\end{equation}
where $a$ is a learnable attention vector and $||$ denotes concatenation.  
Compared to GAT, GATv2 applies the weight matrix before computing the attention score, enabling a more expressive attention mechanism.  
After $T$ layers, we obtain the final node embeddings $h_i^{(T)}$ for all $i \in \mathcal{V}$.  

\subsubsection{Tunable parameter prediction} 
To produce edge-level predictions, node embeddings are combined to form edge embeddings.  
For each transmission line $(i,j) \in \mathcal{E}$, we construct
\begin{equation}
    z_{ij} = [h_i^{(T)} \, || \, h_j^{(T)}],
\end{equation}
which encodes the interaction between its terminal nodes.  
The tunable parameter vector $\tilde{\boldsymbol{\alpha}} \in [0,1]^l$ is then predicted via:
\begin{equation}
    \tilde{\alpha}_{ij} = \text{sigmoid}\left( w_1^\top z_{ij} + w_0 \right),
\end{equation}

\subsection{Parametric Linear Approximation}
Once the tunable parameter $\boldsymbol{\alpha}$ is inferred, the approximate pre-contingency dispatch can be predicted by solving the proposed parametric DC-OPF problem \eqref{eq:p_dcopf}, which can be compactly expressed as
\begin{equation}
\begin{split}
\mathbf{p}_\alpha^\star = \argmin_{\mathbf{p}} \;\; & \frac{1} {2}\mathbf{p^\top C p} + \mathbf{c^\top p} \\
\subjectto \;\; & \mathbf{1^\top p} = b, \; (\lambda) \\
& \mathbf{G p \leq h_\alpha}, \; (\boldsymbol{\mu})
\end{split}
\label{eq:compact_p_dcopf}
\end{equation}
where $\mathbf{p}\in \mathbb{R}^g$ is the dispatch vector, $\lambda\in \mathbb{R}$ and $\boldsymbol{\mu}\in \mathbb{R}^{2(g+l)}$ are the dual variables of the equality and inequality constraints, respectively. $b \in \mathbb{R}$ and $\mathbf{G} \in \mathbb{R}^{2(g+l) \times g}$ are fixed system parameters, and $\mathbf{h}_\alpha \in \mathbb{R}^{2(g+l)}$ encodes the $\alpha$-adjusted line flow constraints. Note that the optimal solution $\mathbf{p}_\alpha^\star$ is a parametric function of the tuning vector $\boldsymbol{\alpha}$ through $\mathbf{h}_\alpha$. For notational convenience, in the remainder of this section, we drop the subscript $\alpha$ that denotes this parameter dependence.

\subsection{Implicit loss function}

The proposed loss function consists of two components that capture both pre- and post-contingency operational costs:
\begin{equation}\label{eq:loss}
    \ell(\boldsymbol{\alpha}) 
    = \ell^{\rm pre}(\boldsymbol{\alpha}) 
    + \ell^{\rm post}\big(\mathbf{p}^\star(\boldsymbol{\alpha})\big),
\end{equation}
where $\ell^{\rm pre}$ represents the optimal value of the parametric linear approximation of the pre-contingency dispatch problem, and $\ell^{\rm post}$ corresponds to the post-contingency cost incurred when this predicted dispatch $\mathbf{p}^\star$ is fixed. Note that $\boldsymbol{\alpha}$ itself depends on the parameters of the GNN. 

To evaluate the quality of a predicted pre-contingency dispatch $\mathbf{p}^\star$, we define an \emph{implicit post-contingency loss function} equal to the optimal objective value of the SC-DCOPF problem.  Formally, the post-contingency loss is given by
\begin{subequations}\label{eq:post_loss}
\begin{align}
\ell^{\rm post}(\mathbf{p}^\star) = \minimize_{\mathbf{p},\mathbf{p}^k,\mathbf{s}^k} \;\;
& \mathbb{E}_\mathbf{d} \Biggl[ \tfrac{\rho}{|\mathcal{K}|} \sum_{k \in \mathcal{K}}\mathbf{1}^\top \mathbf{s}^k \Biggr]\label{eq:loss_a}\\
\rm{s.t.} \;\;
& \forall \mathbf{d}: \nonumber\\
& \qquad \eqref{eq:scopf_e}- \eqref{eq:scopf_i},\; \forall k \in \mathcal{K} \nonumber \\
& \qquad \mathbf{p} = \mathbf{p}^\star,\; (\boldsymbol{\nu}) \label{eq:loss_b}
\end{align}
\end{subequations}
where $\boldsymbol{\nu}$ is the vector of dual variables associated with the fixing constraint \eqref{eq:loss_b}. Note that constraints \eqref{eq:scopf_b}-\eqref{eq:scopf_d} are omitted since $\mathbf{p}^\star$ already satisfies them by construction.

If $\ell^{\rm post}(\mathbf{p}^\star) = 0$, the candidate dispatch $\mathbf{p}^\star$ admits a feasible post-contingency redispatch without requiring load shedding. Conversely, a strictly positive loss indicates the minimum amount of load shedding, averaged across contingencies, that is necessary to restore feasibility when starting from a pre-contingency dispatch $\mathbf{p}^\star$.

\subsection{Backward pass}

During backpropagation, gradients are propagated through the optimization layer by exploiting the differentiability of (i) the implicit loss function \eqref{eq:loss}, and (ii) the KKT conditions of \eqref{eq:compact_p_dcopf}. By the chain rule, the total gradient of the loss is
\begin{equation}
\frac{\partial \ell}{\partial \boldsymbol{\alpha}}
=
\frac{\partial \ell^{\rm pre}}{\partial \boldsymbol{\alpha}}
+
\frac{\partial \ell^{\rm post}}{\partial \mathbf{p}^\star}\,
\frac{\partial \mathbf{p}^\star}{\partial \boldsymbol{\alpha}}.
\end{equation}

\subsubsection{Computing $\frac{\partial \ell^{\rm pre}}{\partial \boldsymbol{\alpha}}$}

The pre-contingency term $\ell^{\rm pre}$ is the optimal value of \eqref{eq:compact_p_dcopf}, where the right-hand side of the inequalities is
$\mathbf{h}$. The gradient of the optimal value with respect to $\mathbf{h}$ is given by the
negative of the optimal dual multipliers on such inequalities:
\begin{equation}
\frac{\partial \ell^{\rm pre}}{\partial \mathbf{h}} \;=\; -\,\boldsymbol{\mu}^\star.
\end{equation}
Therefore, using the chain rule and the fact that $\mathbf{h}$ is a linear
function of $\boldsymbol{\alpha}$,
\begin{equation}
\frac{\partial \ell^{\rm pre}}{\partial \boldsymbol{\alpha}}
=
\bigg(\frac{\partial \mathbf{h}}{\partial \boldsymbol{\alpha}}\bigg)^{\!\top}
\frac{\partial \ell^{\rm pre}}{\partial \mathbf{h}}
\;=\;
-\bigg(\frac{\partial \mathbf{h}}{\partial \boldsymbol{\alpha}}\bigg)^{\!\top}\boldsymbol{\mu}^\star.
\end{equation}

\subsubsection{Computing $\frac{\partial \ell^{\rm post}}{\partial \mathbf{p}^\star}$}

The post-contingency loss function $\ell^{\rm post}(\mathbf{p}^\star)$ is the value function of a convex optimization problem parameterized by $\mathbf{p}^\star$. As such, it inherits differentiability almost everywhere with respect to $\mathbf{p}^\star$. Hence, the dual variables $\boldsymbol{\nu}$ associated with the fixing constraint \eqref{eq:loss_b} can be directly interpreted as the gradient
\begin{equation}
\frac{\partial \ell^{\rm post}}{\partial \mathbf{p}^\star} = -\,\boldsymbol{\nu}^\star,    
\end{equation}
which provides a principled sensitivity of the post-contingency loss with respect to the predicted pre-contingency dispatch. 

\subsubsection{Computing $\frac{\partial \ell^{\rm post}}{\partial \mathbf{p}^\star}\,
\frac{\partial \mathbf{p}^\star}{\partial \boldsymbol{\alpha}}$}

To update the parameters $\boldsymbol{\alpha}$ during backpropagation, we need to compute the gradient of the loss with respect to $\boldsymbol{\alpha}$, which enters the optimization layer through the parameterized bounds $\mathbf{h}$. Thus, the second term of such gradient can be computed by the chain rule as follows:
\begin{equation} \label{eq:gradient_alpha}
\frac{\partial \ell^{\rm post}}{\partial \mathbf{p}^\star}\,
\frac{\partial \mathbf{p}^\star}{\partial \boldsymbol{\alpha}} = 
\frac{\partial \ell^{\rm post}}{\partial \mathbf{p}^\star}\,\frac{\partial \mathbf{p}^\star}{\partial \mathbf{h}} \,
\frac{\partial \mathbf{h}}{\partial \boldsymbol{\alpha}},
\end{equation}
where $\frac{\partial \mathbf{p}^\star}{\partial \mathbf{h}}$ can be computed by differentiating the KKT conditions of \eqref{eq:compact_p_dcopf}. 

The Lagrangian function of \eqref{eq:compact_p_dcopf} is given by
\begin{equation}
\mathcal{L}(\mathbf{p},\lambda,\boldsymbol{\mu})=\frac{1}{2}\mathbf{p^\top Cp}+\mathbf{c^\top p}+\lambda(\mathbf{1^\top p}-b)+\boldsymbol{\mu}^\top(\mathbf{Gp-h})
\end{equation}
where $\lambda$ are the dual variables on the equality constraints
and $\boldsymbol{\mu}\geq \mathbf{0}$ are the dual variables on the inequality constraints.
The KKT conditions of problem \eqref{eq:compact_p_dcopf} are
\begin{equation}
\begin{split}
\mathbf{Cp^\star+c+1\lambda^\star+G^\top}\boldsymbol{\mu}^\star &= \mathbf{0} \\
\mathbf{1^\top p^\star}-b &= 0 \\
\diag(\boldsymbol{\mu}^\star)(\mathbf{Gp^\star-h}) &= \mathbf{0},
\end{split}
\end{equation}
where $\diag(\cdot)$ creates a diagonal matrix from a vector
and $\mathbf{p}^\star$, $\lambda^\star$ and $\boldsymbol{\mu}^\star$ are the optimal
primal and dual variables.
Taking the differentials of these conditions gives the following equations in matrix form
\begin{equation}
  \begin{split}
\boldsymbol{\Gamma}
\begin{bmatrix}
\dd \mathbf{p} \\
\dd \lambda \\
\dd \boldsymbol{\mu} \\
\end{bmatrix} =
-
\begin{bmatrix}
\mathbf{\dd Cp^\star + \dd c + \dd G^\top\boldsymbol{\mu}^\star + \dd 1 \lambda^\star} \\
\mathbf{\dd 1^\top p^\star} - \dd b \\
\mathbf{\diag(\boldsymbol{\mu}^\star) (\dd Gp^\star - \dd h)}
\end{bmatrix},
  \end{split}
  \label{eq:kkt_diff}
\end{equation}
where
\begin{equation}
\boldsymbol{\Gamma} = \begin{bmatrix}
\mathbf{C} & \mathbf{1} & \mathbf{G}^\top  \\
\mathbf{1}^\top & 0 & \mathbf{0}^\top \\
\diag(\boldsymbol{\mu}^\star)\mathbf{G}  & \mathbf{0} & \diag(\mathbf{Gp^\star-h}) \\
\end{bmatrix}
\end{equation}

To compute the sensitivities of the loss function with respect to the tunable parameterized bounds $\mathbf{h}$, we first solve the following adjoint system:
\begin{equation}
\boldsymbol{\Gamma}
\begin{bmatrix}
\mathbf{u}_\mathbf{p} \\
\mathbf{u}_\lambda \\
\mathbf{u}_{\boldsymbol{\mu}} \\
\end{bmatrix} =
\begin{bmatrix}
\left(\frac{\partial \ell^{\rm post}}{\partial \mathbf{p}^\star}\right)^\top \\
0 \\
\mathbf{0}
\end{bmatrix}, \label{eq:adjoint}
\end{equation}
where $\mathbf{u}$ denotes an adjoint vector. Then, the gradient $\frac{\partial \ell^{\rm post}}{\partial \mathbf{p}^\star}\,
\frac{\partial \mathbf{p}^\star}{\partial \mathbf{h}}$ is given by
\begin{equation}
\frac{\partial \ell^{\rm post}}{\partial \mathbf{p}^\star}\,
\frac{\partial \mathbf{p}^\star}{\partial \mathbf{h}} = - \diag(\boldsymbol{\mu}^\star)\mathbf{u}_{\boldsymbol{\mu}}.
\end{equation}

Therefore, using the chain rule and the fact that $\mathbf{h}$ is a linear
function of $\boldsymbol{\alpha}$,
\begin{equation}
\frac{\partial \ell^{\rm post}}{\partial \mathbf{p}^\star}\,
\frac{\partial \mathbf{p}^\star}{\partial \boldsymbol{\alpha}} = -\bigg(\frac{\partial \mathbf{h}}{\partial \boldsymbol{\alpha}}\bigg)^{\!\top}\diag(\boldsymbol{\mu}^\star)\mathbf{u}_{\boldsymbol{\mu}}.
\end{equation}

\section{Numerical Experiments \label{sec:experiments}}

\subsection{Experimental setup}

The performance of the self-supervised model was evaluated through numerical experiments on three benchmark systems: IEEE 57-bus, IEEE 118-bus, and the synthetic Central Illinois 200-bus system. The nominal parameters for these networks were obtained from the Power Grid Library \cite{pglib}. All simulations were conducted on a laptop equipped with an AMD Ryzen AI 9 HX 370 processor (2.00~GHz), Radeon 890M graphics, and 32~GB of RAM.

In the offline stage, a dataset was generated by varying the active power demand within a range of \(\pm 30\% \) with respect to the original load values in the PGLib-OPF instances. For the contingency scenarios, 20\% of the most critical contingencies were considered. The number of variables and constraints of the SC-DCOPF models for all the systems are shown in Table~\ref{tab:model_size_simple}. 
All experiments were implemented in JuMP v1.10, using Gurobi as the optimization solver.

\begin{table}[h]
\centering
\caption{Number of variables and constraints of SC-DCOPF model}
\label{tab:model_size_simple}
\begin{tabular}{lcc}
\hline
\textbf{System} & \textbf{Variables} & \textbf{Constraints} \\ 
\hline
57-bus  & 12,342  & 24,658  \\
118-bus & 61,712  & 127,957 \\
200-bus & 120,900 & 239,829 \\
\hline
\end{tabular}
\end{table}

The GNN was implemented using the GraphNeuralNetworks v0.6.23 package and trained with Flux v0.14.21, while the differentiable optimization layers were implemented in DiffOpt v0.5.0 \cite{Bensacon2024}. The AdamW optimizer was employed with a learning rate of \(1 \times 10^{-6}\). For each test system, a GNN was trained with three graph convolution layers and two hidden layers, using the softplus activation function as a smooth approximation to ReLU. The final layer employed a sigmoid activation function. The same architecture was used across all systems. The network input consisted of node feature vectors of dimension \(2\) and edge feature vectors of dimension \(3\) as detailed in Section III.A. Each attention layer used \(2\) heads and a base hidden dimensionality of \(1024\).

The detailed configuration of the GNN layers is as follows:

\begin{itemize}
    \item \textit{First GATv2Conv layer:} Input dimension (\(2\), \(3\)); output dimension \(1024\) per head. With two heads, the concatenated output has a total dimension of \(2048\).
    \item \textit{Second GATv2Conv layer:} Input dimension (\(2048\), \(3\)); output dimension \(512\) per head. After concatenation, the output dimension is \(1024\).
    \item \textit{Third GATv2Conv layer:} Input dimension (\(1024\), \(3\)); output dimension \(256\) per head. The resulting embeddings have a total dimension of \(512\). After this layer, an \textit{apply\_edges} operation concatenates the features of connected node, resulting in edge representations of size \(1024\).
    \item \textit{Dense layer 1:} Input \(1024\); output \(64\).
    \item \textit{Dense layer 2:} Input \(64\); output \(1\), producing the final prediction.
    \item \textit{Pre-contingency optimization layer:} Solves the parametric DC-OPF for the corresponding operating condition, ensuring that predicted generations satisfy the pre-contingency constraints.
    \item \textit{Post-contingency optimization layer:} Fixes the dispatch decisions after contingency events, enforcing feasibility and system stability under line or component outages.
\end{itemize}


We generated 100 additional samples for model validation, which were not included in the training or testing datasets, to evaluate the proposed model’s out-of-sample performance in terms of cost approximation error, dispatch approximation accuracy, and data efficiency.

\subsection{Numerical Results}

In this section, we compare the performance of the proposed self-supervised framework against three alternative models: semi-supervised, untuned, and end-to-end (E2E) approaches. Each model differs in the level of physical modeling and learning integration, as described below:

\begin{itemize}
    \item \textbf{Self:} The proposed self-supervised framework, where a GNN predicts parameters (\(\alpha\)) that adjust line capacity constraints. The model includes two embedded optimization layers: pre- and post-contingency layers. The loss function is implicitly defined as in \eqref{eq:loss}.
    
    \item \textbf{Semi:} A semi-supervised variant that employs a single optimization layer pertaining to the tunable pre-contingency optimization model \eqref{eq:compact_p_dcopf}. The loss function is the mean-squared error (MSE) between the predicted and true pre-contingency dispatch.
    
    \item \textbf{Untuned:} A baseline that directly uses the dispatches obtained from the untuned DC-OPF formulation, representing a physics-based yet contingency-agnostic model.
    
    \item \textbf{E2E:} An end-to-end GNN model that directly maps input features to pre-contingency nodal power balances through purely data-driven training, without incorporating any optimization layer. The loss function is the MSE between the predicted and true pre-contingency dispatch.
\end{itemize}

\begin{figure*}[ht]
\centering
\subfloat[57-bus system]{\includegraphics[width=0.29\textwidth]{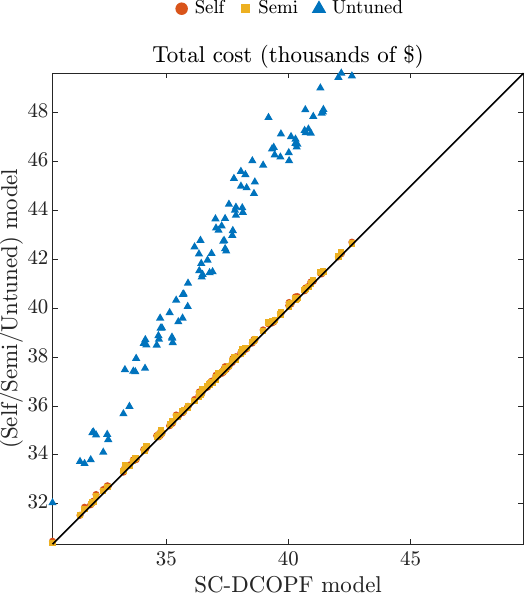}%
\label{fig:cost_57}}
\hfill
\subfloat[118-bus system]{\includegraphics[width=0.29\textwidth]{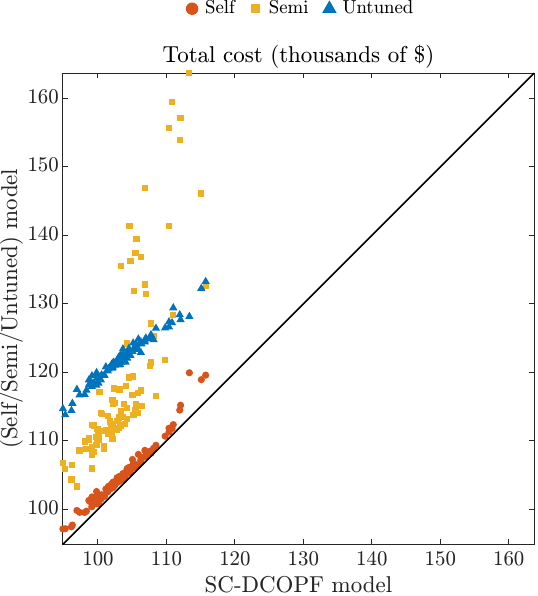}%
\label{fig:cost_118}}
\hfill
\subfloat[200-bus system]{\includegraphics[width=0.29\textwidth]{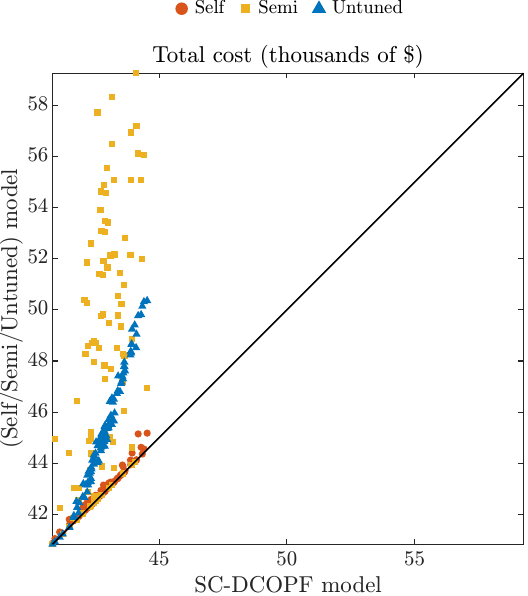}%
\label{fig:cost_200}}
\caption{Cost correlation with respect to SC-DCOPF model.}
\label{fig:cost_correlation}
\end{figure*}

The cost approximation errors, presented in Table~\ref{tab:cost_errors}, are obtained by fixing the pre-contingency dispatch produced by each model (Self, Semi, and Untuned) and solving the post-contingency stage. The self-supervised model achieves the smallest deviations across all test systems, with mean errors below 1.1\%, demonstrating its ability to reproduce near-optimal dispatches even under contingency conditions. In contrast, the semi-supervised model exhibits higher deviations, particularly for the 118-bus and 200-bus systems, indicating sensitivity to the system scale, whereas the untuned baseline consistently shows the largest cost discrepancies. 
Figure~\ref{fig:cost_correlation} illustrates the high cost correlation of the models.

\begin{table}[ht]
\centering
\caption{Cost approximation errors (\%) for each model and test system. Mean (max) values are shown.}
\label{tab:cost_errors}
\renewcommand{\arraystretch}{1.1}
\setlength{\tabcolsep}{6pt}
\begin{tabular}{lccc}
\hline
\multirow{2}{*}{\textbf{Model}} & \multicolumn{3}{c}{\textbf{Test System}} \\
              & \textbf{57-bus} & \textbf{118-bus} & \textbf{200-bus} \\
\hline
Self   & 0.299 (0.789) & \textbf{1.08} (\textbf{5.12}) & \textbf{0.241} (\textbf{2.189}) \\
Semi   & \textbf{0.297} (\textbf{0.770}) & 13.623 (44.054) & 10.512 (31.233) \\
Untuned & 14.192 (22.017) & 17.879 (21.239) & 5.587 (13.331) \\
E2E & 90.141 (286.560) & 201.194 (338.08) & 40.42 (142.83) \\
\hline
\end{tabular}
\end{table}

As shown in Table~\ref{tab:dispatch_accuracy}, the pre-contingency dispatch results across the three test systems indicate that both the self-supervised and semi-supervised models achieve correlations above 0.99 with respect to the true dispatch, demonstrating that the predicted generation profiles closely match the optimal solutions. In particular, the self-supervised model maintains consistently high accuracy across all systems, with mean dispatch errors below 0.07~p.u. even for larger networks. By contrast, the untuned model exhibits noticeably lower correlations and larger deviations due to the lack of information from the contingencies.

\begin{table}[ht]
\centering
\caption{Comparison of pre-contingency dispatch accuracy.}
\renewcommand{\arraystretch}{1.15}
\begin{tabular}{lccc}
\hline
\multirow{2}{*}{\textbf{Model}} & \multirow{2}{*}{\textbf{Correlation}} & \multicolumn{2}{c}{\textbf{Error (p.u.)}} \\ 
& & $\boldsymbol{\rm{mean}}$ & $\boldsymbol{\max}$ \\
\hline
\multicolumn{4}{c}{57-bus} \\
\hline
Self   & 0.9993 & 0.0565 & 0.4311  \\
Semi   & 0.9995 & 0.0477 & 0.3449  \\
Untuned & 0.9143 & 0.8845 & 3.899  \\
\hline
\multicolumn{4}{c}{118-bus} \\
\hline
Self   & 0.9921 & 0.0621 & 3.7195  \\
Semi   & 0.9847 & 0.1198 & 2.9909  \\
Untuned & 0.9683 & 0.1280 & 3.7303 \\
\hline
\multicolumn{4}{c}{200-bus} \\
\hline
Self   & 0.9974 & 0.0192 & 0.4728 \\
Semi   & 0.9986 & 0.0062 & 1.0818  \\
Untuned & 0.9923 & 0.0313 & 1.3457  \\
\hline
\end{tabular}
\label{tab:dispatch_accuracy}
\end{table}

\begin{table}[h]
\centering
\caption{Computation times for dataset generation and model training.}
\renewcommand{\arraystretch}{1.1}
\setlength{\tabcolsep}{4pt}
\begin{tabular}{l | c | c | c c c }
\hline
\textbf{Test} & \textbf{Training} & \textbf{Dataset}&\multicolumn{3}{c}{\textbf{Model Training (min)}} \\
\textbf{System} & \textbf{Samples} & \textbf{(s)} & \textbf{Self} & \textbf{Semi} & \textbf{E2E}  \\
\hline
\multirow{3}{*}{57-bus}       & 1  & 0.159 & 1.245 & 0.948 & 1.0944 \\
             & 10 & 1.771 & 4.501 & 7.076 & 2.539  \\
             & 25 & 4.199 & 20.608    & 11.572    & 8.057     \\
             & 100 & 15.162 & --    & --    & 25.194     \\
\hline
\multirow{3}{*}{118-bus} & 1  & 1.261 & 12.728 & 16.238 & 0.776 \\
             & 10 & 14.168 & 291.580 & 234.650 & 4.541  \\
             & 25 & 34.768 & 334.430 & 299.980 & 11.644  \\
             & 100 & 132.56 & --    & --    & 39.624     \\
\hline
\multirow{3}{*}{200-bus} & 1  & 2.993 & 11.839 & 22.967 & 0.752  \\
             & 10 & 28.397 & 132.920 & 89.330 & 6.833  \\
             & 25 & 68.395 & 272.680 & 115.550 & 15.871 \\
             & 100 & 276.92 & --    & --    & 49.902     \\
\hline
\end{tabular}
\label{tab:training_times}
\end{table}

Table~\ref{tab:training_times} reports the offline computational time for all learning-based models. As expected, both dataset generation and training times increase with the system size and the number of training samples. The self- and semi-supervised frameworks require longer training due to the embedded optimization layers, whereas the end-to-end (E2E) model trains considerably faster. 


\subsection{Data Efficiency Comparison}

Table~\ref{tab:objective_errors} summarizes the mean and maximum relative errors of the objective function for different numbers of training samples. The self-supervised model shows a clear advantage over the other models, achieving high accuracy even when trained with very limited training data. For instance, in the 118-bus and 200-bus systems, the self-supervised model reaches mean cost errors below 2\% with only 10 training samples, while the semi-supervised and E2E models require significantly more data to achieve comparable accuracy. 

 These results demonstrate the data efficiency of the proposed self-supervised learning framework, which leverages the underlying physical structure of the pre- and post-contingency stages through the optimization layers. In contrast, the semi-supervised and end-to-end models rely solely on pre-contingency physics and/or labeled data, resulting in slower convergence and higher sensitivity to sample size. In general, the proposed self-supervised approach not only reduces data requirements but also delivers consistently better performance across different test systems.

\begin{table}[ht]
\centering
\caption{Mean (max) relative errors [\%] of the objective function}
\renewcommand{\arraystretch}{1.1}
\setlength{\tabcolsep}{5pt}
\begin{tabular}{l|c|ccc}
\hline
\textbf{Test} & \textbf{Training} & \multicolumn{3}{|c}{\textbf{Model}} \\
\textbf{System} & \textbf{Samples} & \textbf{Self} & \textbf{Semi} & \textbf{E2E} \\
\hline
\multirow{4}{*}{57-bus} 
        & 1   & 0.49 (1.11) & \textbf{0.43} (\textbf{1.02}) & 98.98 (243.46) \\
        & 10  & 0.30 (0.79) & \textbf{0.30} \textbf{(0.77)} & 89.43 (322.56) \\
        & 25  &  0.30 (0.78) & \textbf{0.29} \textbf{(0.75)} & 90.14 (286.56) \\
        & 100 & -- & -- & 92.38 (267.64) \\
\hline
\multirow{4}{*}{118-bus} 
        & 1   & \textbf{15.65} (\textbf{44.92}) & 50.44 (94.03) & 293.90 (441.82) \\
        & 10  & \textbf{1.40} (\textbf{5.81})   & 14.93 (44.41) & 270.73 (417.02) \\
        & 25  & \textbf{1.08} (\textbf{5.12})   & 13.62 (44.05) & 201.20 (338.09) \\
        & 100 & --            & --            & 376.21 (533.67) \\
\hline
\multirow{4}{*}{200-bus} 
        & 1   & \textbf{16.57} (\textbf{42.35}) & 17.96 (44.16) & 44.60 (153.37) \\
        & 10  & \textbf{0.24} (\textbf{2.19})   & 12.66 (35.58) & 41.18 (144.76) \\
        & 25  & \textbf{0.24} (\textbf{2.44})   & 10.51 (31.23) & 40.42 (142.81) \\
        & 100 & --            & --            & 40.42 (142.83) \\
\hline
\end{tabular}
\label{tab:objective_errors}
\end{table}

\section{Conclusions and Future Work\label{sec:conclusion}}

This paper introduced a self-supervised learning framework for approximating the Security-Constrained DC Optimal Power Flow (SC-DCOPF) problem using a tunable parametric linear model. By leveraging demand-dependent parameters and embedding differentiable optimization layers, the proposed method achieves high dispatch accuracy and cost efficiency without requiring labeled SC-DCOPF solutions. Numerical experiments across multiple benchmark systems demonstrate the model’s scalability, data efficiency, and better performance compared to semi-supervised and end-to-end baselines.

Future work will focus on enhancing scalability by exploiting the decomposability of the implicit loss function, enabling parallel training across contingencies using GPU-based solvers. The tuned parametric models can also be integrated with decomposition techniques such as Benders decomposition to accelerate convergence by improving solution quality during early iterations. Additionally, extending this framework to incorporate AC power flow models, and dynamic system behaviors could enhance its applicability to real-world operations. Another promising direction is the development of mixed-integer parametric models to capture discrete operational decisions such as unit commitment and topology control. Investigating robustness under uncertainty, and integrating real-time data streams and forecast errors, will further support secure and adaptive power system operations.

\bibliography{bib/references.bib}
\bibliographystyle{IEEEtran}


\end{document}